\documentclass[preprint,12pt]{elsarticle}
\usepackage[figuresright]{rotating}
\usepackage{amssymb}
\usepackage{amsthm}
\usepackage{multicol}
\usepackage{subfigure}
\usepackage{graphicx}
\usepackage{epstopdf}
\usepackage{fullpage}
\usepackage{latexsym,amsmath}
\usepackage{stmaryrd}
\usepackage{algorithm,algorithmic}
\usepackage{subfigure}
\usepackage{amsfonts}
\usepackage{xcolor,multirow}
\usepackage{bm}
\usepackage{url}
\usepackage{diagbox}
\usepackage{array}
\usepackage{booktabs}
\usepackage{wrapfig} 
\biboptions{sort&compress}

\newtheorem{theorem}{Theorem}
\newtheorem{definition}{Definition}
\newtheorem{proposition}{Proposition}

\usepackage{amssymb}

\usepackage{tikz}  
\usetikzlibrary{arrows.meta}

\begin{document}
\begin{sloppypar}
\begin{frontmatter}



\title{Quaternion tensor singular value decomposition using a flexible transform-based approach}

\author[lab1]{Jifei Miao}
\ead{jifmiao@163.com}
\author[lab1]{Kit Ian Kou\corref{cor1}}
\ead{kikou@umac.mo}

\address[lab1]{Department of Mathematics, Faculty of
	Science and Technology, University of Macau, Macau 999078, China}
\cortext[cor1]{Corresponding author}

\begin{abstract}
A flexible transform-based tensor product named $\star_{{\rm{QT}}}$-product for $L$th-order ($L\geq 3$) quaternion tensors is proposed. Based on the $\star_{{\rm{QT}}}$-product, we define the corresponding singular value decomposition named TQt-SVD and the rank named TQt-rank of the $L$th-order ($L\geq 3$) quaternion tensor. Furthermore, with orthogonal quaternion transformations, the TQt-SVD can provide the best TQt-rank-$s$ approximation of any $L$th-order ($L\geq 3$) quaternion tensor. In the experiments, we have verified the effectiveness of the proposed TQt-SVD in the application of the best TQt-rank-$s$ approximation for color videos represented by third-order quaternion tensors.
\end{abstract}

%


\end{frontmatter}

\section{Introduction}
As the quaternion can well preserve cross-channel correlation of color channels, quaternion matrices and quaternion tensors have recently been widely used in the processing of color images and color videos \cite{DBLP:journals/tip/ChenXZ20,DBLP:journals/nla/JiaNS19,DBLP:journals/pr/MiaoKL20}, which have proven that the quaternion is a very suitable tool for color pixel representation. To make low-rank approximations for color videos, the authors in \cite{DBLP:journals/appml/QinMZ22}  generalized the tensor singular value decomposition (t-SVD) \cite{KILMER2011641} to the quaternion domain and proposed a quaternion tensor singular value decomposition (Qt-SVD).
Following the important property, the same as t-product \cite{KILMER2011641} based on, that any complex circulant matrix can be diagonalized by the normalized Discrete Fourier Transformation (DFT) matrix, the authors in \cite{DBLP:journals/appml/QinMZ22} designed a Qt-product by rewriting the quaternion tensor as two complex tensors. However, the Qt-SVD is only applicable to third-order
quaternion tensors. In addition, the transform that the Qt-product depends on can only be fixed as DFT, which may cause inflexible applications and the inability to adapt to different color video data. 

In this paper, to handle the limitations of Qt-SVD, inspired by the real-valued cases in \cite{KERNFELD2015545,Kilmere2015851118}, we first introduce a product for any $L$th-order ($L\geq 3$) quaternion tensors ($\star_{{\rm{QT}}}$-product) based on any invertible quaternion transform. Then, following the $\star_{{\rm{QT}}}$-product, we propose a flexible transform-based Qt-SVD named TQt-SVD. The ``flexible'' of  TQt-SVD means that it is applicable for any $L$th-order ($L\geq 3$) quaternion tensor (not necessary third-order) by using any invertible quaternion transform (not necessary fixed one). Furthermore, with  orthogonal quaternion transformations, the TQt-SVD can provide the best TQt-rank-$s$ approximation of a quaternion tensor, which provides a new way to process any higher-dimensional visual data with three color channels, \emph{e.g.}, color videos. Since the orthogonal quaternion transformation does not need to be fixed, we can select different transformations for each color video, or use data-driven transformations. The experiments verify the effectiveness of the proposed TQt-SVD in the application of the best TQt-rank-$s$ approximation for color videos.

\section{Quaternion-related knowledge and some notations}
Some quaternion-related knowledge including quaternion matrices and tensors can be found in the supplementary material. In this paper, $\mathbb{R}$, $\mathbb{C}$, and $\mathbb{H}$ respectively denote the real, complex, and quaternion spaces. A quaternion scalar, vector, matrix, and tensor are written as $\dot{a}$, $\dot{\mathbf{a}}$, $\dot{\mathbf{A}}$, and $\dot{\mathcal{A}}$, respectively. For an $L$th-order ($L\geq 3$) quaternion tensor $\dot{\mathcal{A}}\in\mathbb{H}^{N_{1}\times N_{2} \times N_{3}\times\ldots \times N_{L}}$, $\dot{\mathcal{A}}(:,:,n_{3},\ldots,n_{L})$ (where $n_{k}=1, \ldots, N_{k}$, for $k=3, \ldots, L$)  denotes its frontal slices, and $\dot{\mathbf{A}}_{(k)}$ denotes its mode-$k$ unfolding. The $\times_{k}$ denotes the $k$-mode product, $[K]:=\{1,2,\ldots,K\}$, and $\|\cdot\|_{F}$ is the Frobenius norm.
The f-diagonal quaternion tensor means that its
entries only lie along the diagonal of its frontal slices.

\section{The $\star_{{\rm{QT}}}$-product and TQt-SVD}
\label{sec:2}
In this section, we first introduce the $\star_{{\rm{QT}}}$-product and the associated concepts,  then propose the TQt-SVD.
\begin{definition}(Quaternion facewise product ($\star_{{\rm{QF}}}$-product)):
	Given two $L$th-order ($L\geq 3$) quaternion tensors $\dot{\mathcal{A}}\in\mathbb{H}^{N_{1}\times P \times N_{3}\times\ldots \times N_{L}}$ and $\dot{\mathcal{B}}\in\mathbb{H}^{P\times N_{2}\times N_{3} \times\ldots \times N_{L}}$, the $\star_{{\rm{QF}}}$-product can be written as
\begin{equation}\label{equ1}
	\dot{\mathcal{F}}=\dot{\mathcal{A}}\star_{{\rm{QF}}}\dot{\mathcal{B}}\in\mathbb{H}^{N_{1}\times N_{2} \times\ldots \times N_{L}},
\end{equation}
such that each frontal slice of $\dot{\mathcal{F}}$ satisfies
$\dot{\mathcal{F}}(:,:,n_{3},\ldots,n_{L})=	\dot{\mathcal{A}}(:,:,n_{3},\ldots,n_{L})	\dot{\mathcal{B}}(:,:,n_{3},\ldots,n_{L})$, where $n_{k}=1, \ldots, N_{k}$, for $k=3, \ldots, L$.
\end{definition}
Based on the above definition of $\star_{{\rm{QF}}}$-product, we define the following $\star_{{\rm{QT}}}$-product.
\begin{definition}($\star_{{\rm{QT}}}$-product):
	Given two $L$th-order ($L\geq 3$) quaternion tensors $\dot{\mathcal{A}}\in\mathbb{H}^{N_{1}\times P \times N_{3}\times\ldots \times N_{L}}$ and $\dot{\mathcal{B}}\in\mathbb{H}^{P\times N_{2}\times N_{3} \times\ldots \times N_{L}}$, with $L-2$ invertible quaternion matrices $\dot{\mathbf{T}}_{3}\in\mathbb{H}^{N_{3}\times N_{3}},\ldots, \dot{\mathbf{T}}_{L}\in\mathbb{H}^{N_{L}\times N_{L}}$,
	the $\star_{{\rm{QT}}}$-product can be written as
	\begin{equation}\label{equ2}
		\dot{\mathcal{T}}=\dot{\mathcal{A}}\star_{{\rm{QT}}}\dot{\mathcal{B}}=(\widehat{\dot{\mathcal{A}}}\star_{{\rm{QF}}}\widehat{\dot{\mathcal{B}}})\times_{3}\dot{\mathbf{T}}_{3}^{-1}\times\cdots\times_{L}\dot{\mathbf{T}}_{L}^{-1},
	\end{equation}
	where $\widehat{\dot{\mathcal{A}}}=\dot{\mathcal{A}}\times_{3}\dot{\mathbf{T}}_{3}\times\cdots\times_{L}\dot{\mathbf{T}_{L}}$ and $\widehat{\dot{\mathcal{B}}}=\dot{\mathcal{B}}\times_{3}\dot{\mathbf{T}}_{3}\times\cdots\times_{L}\dot{\mathbf{T}}_{L}$.
\end{definition}

We then define the conjugate transpose, the identity quaternion tensor, and the unitary quaternion tensor under the above defined $\star_{{\rm{QT}}}$-product.

\begin{definition}(Conjugate transpose):
Given an $L$th-order ($L\geq 3$) quaternion tensor $\dot{\mathcal{A}}\in\mathbb{H}^{N_{1}\times N_{2} \times\ldots \times N_{L}}$, the conjugate transpose is 
denoted as $\dot{\mathcal{A}}^{H}\in\mathbb{H}^{N_{2}\times N_{1} \times\ldots \times N_{L}}$, satisfying
\begin{equation*}
\widehat{\dot{\mathcal{A}}}^{H}(:,:,n_{3},\ldots,n_{L})	=(\widehat{\dot{\mathcal{A}}}(:,:,n_{3},\ldots,n_{L}))^{H},
\end{equation*} 
where $n_{k}=1, \ldots, N_{k}$, for $k=3, \ldots, L$.
\end{definition}

\begin{definition}(Identity quaternion tensor):\label{iqt}
	 The identity $L$th-order ($L\geq 3$) quaternion tensor $\dot{\mathcal{I}}\in\mathbb{H}^{P\times P\times N_{3} \times\ldots \times N_{L}}$ is the quaternion tensor satisfying that each frontal
	 slice of $\widehat{\dot{\mathcal{I}}}$ is the identity quaternion matrix, i.e.,
	 \begin{equation}\label{equ3}
	 \widehat{\dot{\mathcal{I}}}(:,:,n_{3},\ldots,n_{L})=\dot{\mathbf{I}}\in\mathbb{H}^{P\times P},
	 \end{equation} 
 where $n_{k}=1, \ldots, N_{k}$, for $k=3, \ldots, L$.
\end{definition}

\begin{proposition}
	Let $\dot{\mathcal{I}}\in\mathbb{H}^{P\times P\times N_{3} \times\ldots \times N_{L}}$ be the identity $L$th-order ($L\geq 3$) quaternion tensor defined in Definition \ref{iqt}. Then, for any $L$th-order quaternion tensors $\dot{\mathcal{A}}\in\mathbb{H}^{P\times N_{2}\times N_{3} \times\ldots \times N_{L}}$ and $\dot{\mathcal{B}}\in\mathbb{H}^{N_{1}\times P\times N_{3} \times\ldots \times N_{L}}$, we have
	\begin{equation*}
		\dot{\mathcal{I}}\star_{{\rm{QT}}}\dot{\mathcal{A}}=\dot{\mathcal{A}}\qquad \text{and}\qquad \dot{\mathcal{B}}\star_{{\rm{QT}}}\dot{\mathcal{I}}=\dot{\mathcal{B}}.
	\end{equation*}
\end{proposition}
\textit{Proof.} Denote $\dot{\mathcal{Y}}=\widehat{\dot{\mathcal{I}}}\star_{{\rm{QF}}}\widehat{\dot{\mathcal{A}}}$, then $\dot{\mathcal{Y}}(:,:,n_{3},\ldots,n_{L})=\widehat{\dot{\mathcal{I}}}(:,:,n_{3},\ldots,n_{L})	\widehat{\dot{\mathcal{A}}}(:,:,n_{3},\ldots,n_{L})=\dot{\mathbf{I}}\widehat{\dot{\mathcal{A}}}(:,:,n_{3},\ldots,n_{L})=\widehat{\dot{\mathcal{A}}}(:,:,n_{3},\ldots,n_{L})$ for $n_{k}=1, \ldots, N_{k}$, and $k=3, \ldots, L$, \emph{i.e.}, $\widehat{\dot{\mathcal{I}}}\star_{{\rm{QF}}}\widehat{\dot{\mathcal{A}}}=\widehat{\dot{\mathcal{A}}}$. From (\ref{equ2}), we have $\widehat{(\dot{\mathcal{I}}\star_{{\rm{QT}}}\dot{\mathcal{A}})}=\widehat{\dot{\mathcal{I}}}\star_{{\rm{QF}}}\widehat{\dot{\mathcal{A}}}=\widehat{\dot{\mathcal{A}}}$, which illustrates that $\dot{\mathcal{I}}\star_{{\rm{QT}}}\dot{\mathcal{A}}=\dot{\mathcal{A}}$. Similarly, we can get $\dot{\mathcal{B}}\star_{{\rm{QT}}}\dot{\mathcal{I}}=\dot{\mathcal{B}}$.

\begin{definition}(Unitary
	quaternion tensor):
	An $L$th-order ($L\geq 3$) quaternion tensor $\dot{\mathcal{U}}\in\mathbb{H}^{P\times P\times N_{3} \times\ldots \times N_{L}}$ is unitary if
	\begin{equation*}
		\dot{\mathcal{U}}^{H}\star_{{\rm{QT}}}\dot{\mathcal{U}}=\dot{\mathcal{I}}=\dot{\mathcal{U}}\star_{{\rm{QT}}}\dot{\mathcal{U}}^{H}.
	\end{equation*}
\end{definition}

Now we propose the following TQt-SVD for any $L$th-order ($L\geq 3$) quaternion tensor.
\begin{theorem}(TQt-SVD): Any $L$th-order ($L\geq 3$) quaternion tensor $\dot{\mathcal{A}}\in\mathbb{H}^{N_{1}\times N_{2} \times\ldots \times N_{L}}$ 
can be decomposed as
\begin{equation}\label{equ4}
	\dot{\mathcal{A}}=\dot{\mathcal{U}}\star_{{\rm{QT}}}\dot{\mathcal{D}}\star_{{\rm{QT}}} \dot{\mathcal{V}}^{H},
\end{equation}
 where $\dot{\mathcal{U}}\in\mathbb{H}^{N_{1}\times N_{1}\times N_{3} \times\ldots \times N_{L}}$ and $\dot{\mathcal{V}}\in\mathbb{H}^{N_{2}\times N_{2}\times N_{3} \times\ldots \times N_{L}}$	are unitary 
 quaternion tensors, and $\dot{\mathcal{D}}\in\mathbb{H}^{N_{1}\times N_{2}\times N_{3} \times\ldots \times N_{L}}$ is an f-diagonal quaternion tensor.
\end{theorem}
\textit{Proof.} Based on the quaternion matrix SVD \cite{1997Quaternions},  we can decompose $\widehat{\dot{\mathcal{A}}}(:,:,n_{3},\ldots,n_{L})$ by
\begin{equation}\label{equ5}
\widehat{\dot{\mathcal{A}}}(:,:,n_{3},\ldots,n_{L})=\widehat{\dot{\mathcal{U}}}(:,:,n_{3},\ldots,n_{L})\widehat{\dot{\mathcal{D}}}(:,:,n_{3},\ldots,n_{L})\widehat{\dot{\mathcal{V}}}(:,:,n_{3},\ldots,n_{L})^{H},	
\end{equation}
where $\widehat{\dot{\mathcal{U}}}(:,:,n_{3},\ldots,n_{L})$ and  $\widehat{\dot{\mathcal{V}}}(:,:,n_{3},\ldots,n_{L})$ are unitary quaternion matrices, $\widehat{\dot{\mathcal{D}}}(:,:,n_{3},\ldots,n_{L})$ is diagonal
for all $n_{k}=1, \ldots, N_{k}$, and $k=3, \ldots, L$. Denote $\dot{\mathcal{X}}=\widehat{\dot{\mathcal{U}}}^{H}\star_{{\rm{QF}}}\widehat{\dot{\mathcal{U}}}$, then $\dot{\mathcal{X}}(:,:,n_{3},\ldots,n_{L})=\widehat{\dot{\mathcal{U}}}(:,:,n_{3},\ldots,n_{L})^{H}\widehat{\dot{\mathcal{U}}}(:,:,n_{3},\ldots,n_{L})=\dot{\mathbf{I}}=\widehat{\dot{\mathcal{I}}}(:,:,n_{3},\ldots,n_{L})$ for all $n_{k}=1, \ldots, N_{k}$, and $k=3, \ldots, L$, which means $\widehat{\dot{\mathcal{X}}}=\dot{\mathcal{I}}$. Hence, $\widehat{(\dot{\mathcal{U}}^{H}\star_{{\rm{QT}}}\dot{\mathcal{U}})}=\widehat{\dot{\mathcal{U}}}^{H}\star_{{\rm{QF}}}\widehat{\dot{\mathcal{U}}}=\dot{\mathcal{X}}=\widehat{\dot{\mathcal{I}}}$, which illustrates that $\dot{\mathcal{U}}^{H}\star_{{\rm{QT}}}\dot{\mathcal{U}}=\dot{\mathcal{I}}$. Similarly, we can get $\dot{\mathcal{U}}\star_{{\rm{QT}}}\dot{\mathcal{U}}^{H}=\dot{\mathcal{I}}$, and $\dot{\mathcal{V}}^{H}\star_{{\rm{QT}}}\dot{\mathcal{V}}=\dot{\mathcal{I}}=\dot{\mathcal{V}}\star_{{\rm{QT}}}\dot{\mathcal{V}}^{H}$. Furthermore, it is clear that $\dot{\mathcal{D}}$ is  f-diagonal, since  $\widehat{\dot{\mathcal{D}}}(:,:,n_{3},\ldots,n_{L})$ is diagonal
for all $n_{k}=1, \ldots, N_{k}$, and $k=3, \ldots, L$, and note that for any invertible quaternion matrices $\dot{\mathbf{T}}_{3}\in\mathbb{H}^{N_{3}\times N_{3}},\ldots, \dot{\mathbf{T}}_{L}\in\mathbb{H}^{N_{L}\times N_{L}}$, $0\times_{3}\dot{\mathbf{T}}_{3}\times\cdots\times_{L}\dot{\mathbf{T}_{L}}=0=0\times_{3}\dot{\mathbf{T}}_{3}^{-1}\times\cdots\times_{L}\dot{\mathbf{T}_{L}}^{-1}$. 

Based on the TQt-SVD, we define the following rank for $L$th-order ($L\geq 3$) quaternion tensor.
\begin{definition}(TQt-rank):
 For any $L$th-order ($L\geq 3$) quaternion tensor $\dot{\mathcal{A}}\in\mathbb{H}^{N_{1}\times N_{2} \times\ldots \times N_{L}}$ and its TQt-SVD $\dot{\mathcal{A}}=\dot{\mathcal{U}}\star_{{\rm{QT}}}\dot{\mathcal{D}}\star_{{\rm{QT}}} \dot{\mathcal{V}}^{H}$, the TQt-rank of $\dot{\mathcal{A}}$ is the number of nonzero tubes in $\dot{\mathcal{D}}$, i.e.,
 \begin{equation*}
 	rank_{TQt}(\dot{\mathcal{A}})=\#\{k| \|\dot{\mathcal{D}}(k,k,:,\ldots,:)\|_{F}>0\},\ k\in[K],\ K=min(N_{1},N_{2}).
 \end{equation*}
Furthermore, the $k$th singular value of $\dot{\mathcal{A}}$ is defined as $\sigma_{k}=\|\dot{\mathcal{D}}(k,k,:,\ldots,:)\|_{F}$ for $k\in[K]$.
\end{definition}

For any $L$th-order ($L\geq 3$) quaternion tensor $\dot{\mathcal{A}}\in\mathbb{H}^{N_{1}\times N_{2} \times\ldots \times N_{L}}$ and its TQt-SVD $\dot{\mathcal{A}}=\dot{\mathcal{U}}\star_{{\rm{QT}}}\dot{\mathcal{D}}\star_{{\rm{QT}}} \dot{\mathcal{V}}^{H}$,
supposing that $rank_{TQt}(\dot{\mathcal{A}})=r$, then $\dot{\mathcal{A}}$ can be rewritten as
\begin{equation*}
	\dot{\mathcal{A}}=\sum_{k=1}^{r}\dot{\mathcal{U}}(:,k,:,\ldots,:)\star_{{\rm{QT}}}\dot{\mathcal{D}}(k,k,:,\ldots,:)\star_{{\rm{QT}}} \dot{\mathcal{V}}(:,k,:,\ldots,:)^{H}.
\end{equation*}

\begin{theorem}\label{them2} Let $\dot{\mathcal{A}}\in\mathbb{H}^{N_{1}\times N_{2} \times\ldots \times N_{L}}$ be any $L$th-order ($L\geq 3$) quaternion tensor with the TQt-SVD $\dot{\mathcal{A}}=\dot{\mathcal{U}}\star_{{\rm{QT}}}\dot{\mathcal{D}}\star_{{\rm{QT}}} \dot{\mathcal{V}}^{H}$ where the $\star_{{\rm{QT}}}$-product consists of only multiples of orthogonal quaternion transformations, i.e., $\dot{\mathbf{T}}_{n}=c_{n}\dot{\mathbf{U}}_{n}\in\mathbb{H}^{N_{n}\times N_{n}}, n=3,\ldots,L$, where $c_{n}\in\mathbb{R}\neq 0$, $\dot{\mathbf{U}}_{n}\in\mathbb{H}^{N_{n}\times N_{n}}$ is unitary for $n=3,\ldots,L$. Define $\dot{\mathcal{A}}_{s}=\sum_{k=1}^{s}\dot{\mathcal{U}}(:,k,:,\ldots,:)\star_{{\rm{QT}}}\dot{\mathcal{D}}(k,k,:,\ldots,:)\star_{{\rm{QT}}} \dot{\mathcal{V}}(:,k,:,\ldots,:)^{H}$ as the TQt-rank-$s$ approximation of $\dot{\mathcal{A}}$. Then,
\begin{equation*}
\dot{\mathcal{A}}_{s}=\mathop{{\rm{arg\, min}}}\limits_{\dot{\mathcal{G}}\in\mathbb{H}^{N_{1}\times N_{2} \times\ldots \times N_{L}}}\ \|\dot{\mathcal{A}}-\dot{\mathcal{G}}\|_{F}, \quad \text{where} \quad rank_{TQt}(\dot{\mathcal{G}})=s,
\end{equation*}
i.e., $\dot{\mathcal{A}}_{s}$ is the best TQt-rank-$s$ approximation of $\dot{\mathcal{A}}$. Moreover, the squared error between $\dot{\mathcal{A}}$ and $\dot{\mathcal{A}}_{s}$ is $\|\dot{\mathcal{A}}-\dot{\mathcal{A}}_{s}\|_{F}^{2}=\sum_{k=s+1}^{r}\|\dot{\mathcal{D}}(k,k,:,\ldots,:)\|_{F}^{2}$, where $r$ is the TQt-rank of $\dot{\mathcal{A}}$.
\end{theorem}

The proof of Theorem \ref{them2} can be found in the supplementary material.

\section{Experiments}
In the experiments, we aim to verify the effectiveness of the proposed TQt-SVD in the application of the best TQt-rank-$s$ approximation for color videos. A color video with RGB three channels is represented as a pure third-order quaternion tensor $\dot{\mathcal{Q}}=\mathcal{Q}_{1}i+\mathcal{Q}_{2}j+\mathcal{Q}_{3}k\in\mathbb{H}^{N_{1}\times N_{2} \times F}$, where the size of each frame of the color video is $N_{1}\times N_{2}$, and $F$ denotes the number of frames. Specifically, the RGB three channels of the pixel values are set as the $i$, $j$, and $k$ parts of the quaternion tensor $\dot{\mathcal{Q}}$, respectively. Thus, the real part of $\dot{\mathcal{Q}}$ vanishes. Let the TQt-SVD of $\dot{\mathcal{Q}}$ be $\dot{\mathcal{Q}}=\dot{\mathcal{U}}\star_{{\rm{QT}}}\dot{\mathcal{D}}\star_{{\rm{QT}}} \dot{\mathcal{V}}^{H}$ where the $\star_{{\rm{QT}}}$-product consists of an orthogonal quaternion transformation $\dot{\mathbf{T}}=c\dot{\mathbf{U}}\in\mathbb{H}^{F\times F}$, where $c\in\mathbb{R}\neq 0$, $\dot{\mathbf{U}}\in\mathbb{H}^{F\times F}$ is unitary.
Then, from Theorem \ref{them2}, the best TQt-rank-$s$ approximation for the color video is $
\dot{\mathcal{Q}}_{s}=\sum_{k=1}^{s}\dot{\mathcal{U}}(:,k,:)\star_{{\rm{QT}}}\dot{\mathcal{D}}(k,k,:)\star_{{\rm{QT}}} \dot{\mathcal{V}}(:,k,:)^{H}$.	

\begin{figure}[htbp]
	\centering
	\subfigure[AN119T]{\includegraphics[width=2.5cm,height=1.3cm]{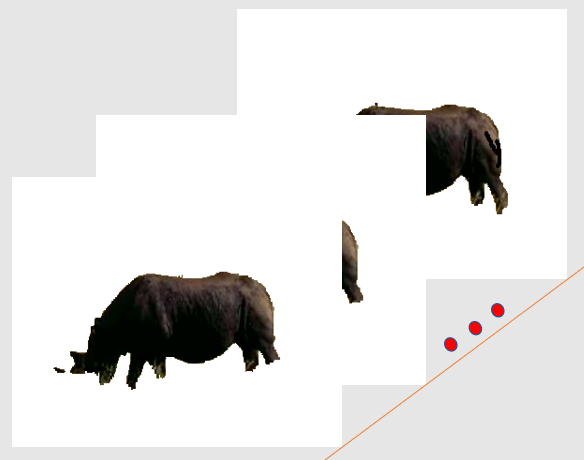}}\hspace{1cm}
	\subfigure[DO01$\_$013]{\includegraphics[width=2.5cm,height=1.3cm]{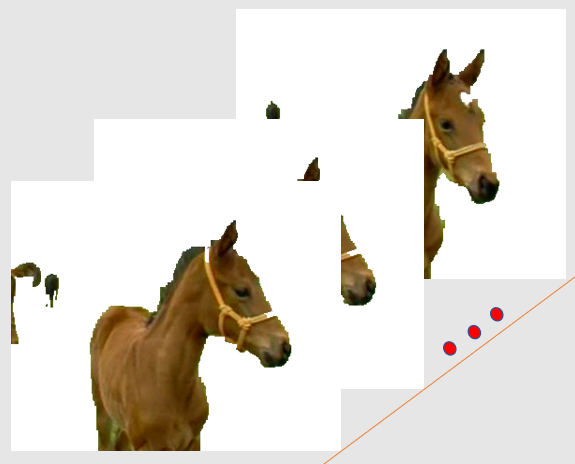}}\hspace{1cm}
	\subfigure[DO01$\_$030]{\includegraphics[width=2.5cm,height=1.3cm]{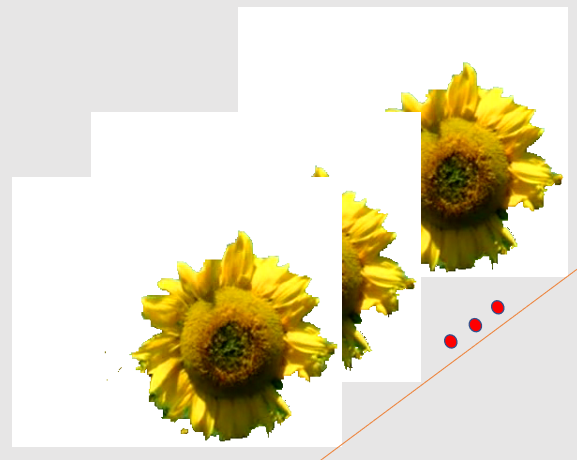}}\hspace{1cm}
	\subfigure[DO01$\_$014]{\includegraphics[width=2.5cm,height=1.3cm]{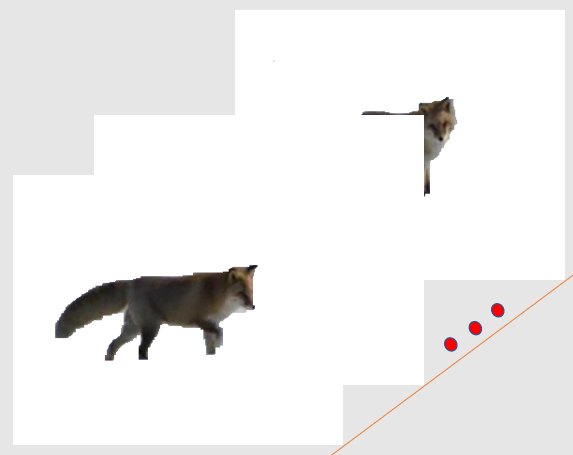}}
	\caption{The $4$ tested color videos, reconstructed by third-order quaternion tensors.}
	\label{fig_1}
\end{figure}

\begin{figure}[htbp]
	\centering
    \subfigure[The rank-$20$ approximation of one randomly selected frame of the AN119T.]{\includegraphics[width=14.5cm,height=2.3cm]{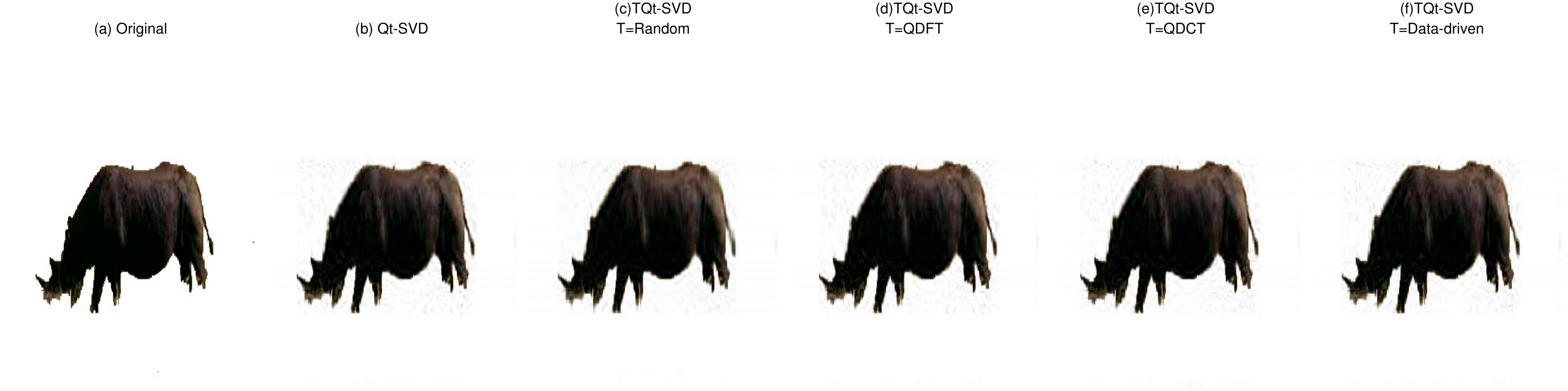}}
	\subfigure[The rank-$20$ approximation of one randomly selected frame of the DO01$\_$013.]{\includegraphics[width=14.5cm,height=2.3cm]{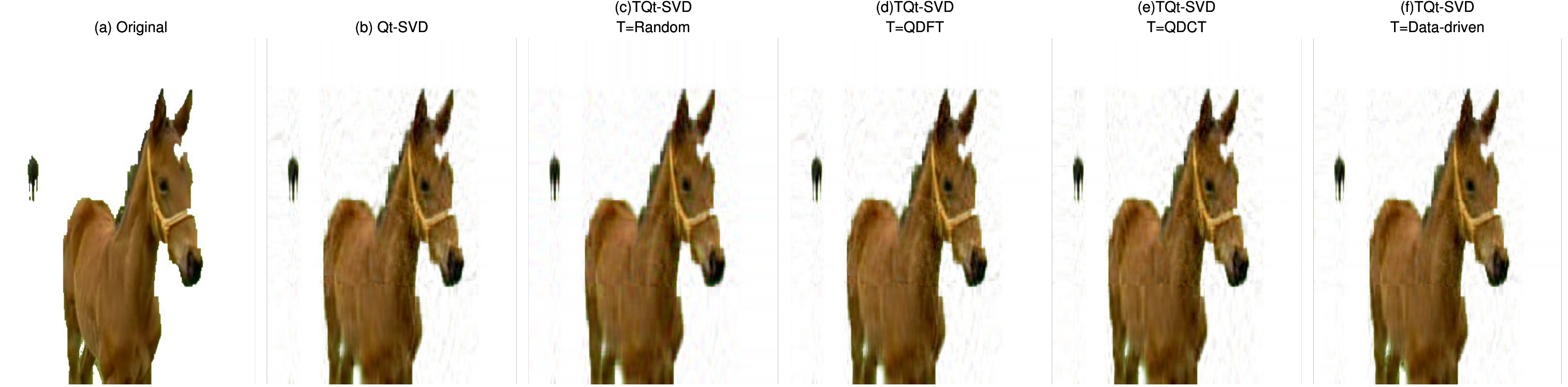}}
	\subfigure[The rank-$20$ approximation of one randomly selected frame of the DO01$\_$014.]{\includegraphics[width=14.5cm,height=2.3cm]{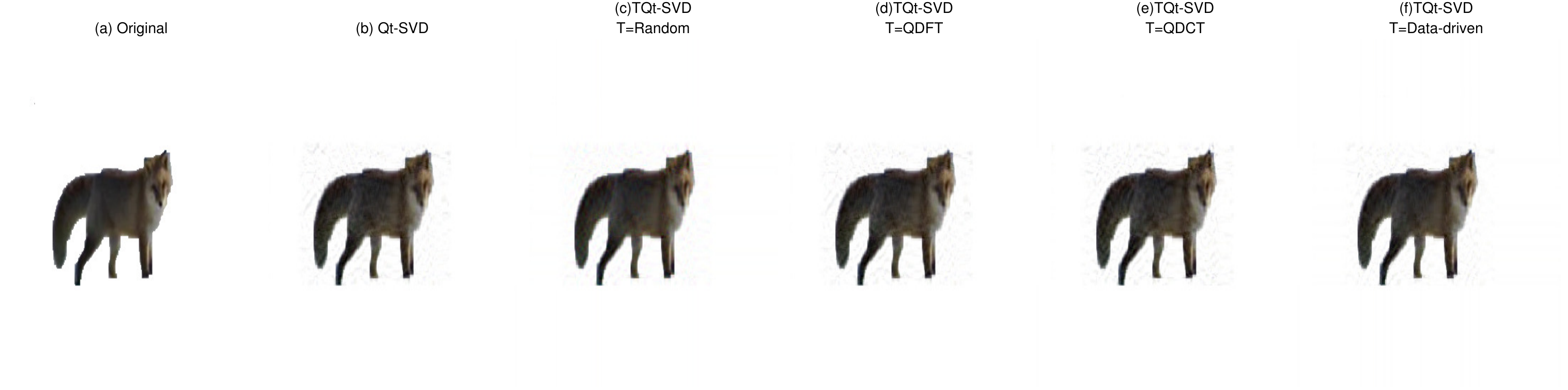}}
	\caption{The rank-20 approximation of one randomly selected frame of the (i) AN119T, (ii) DO01$\_$013, and (iii) DO01$\_$014.}
	\label{fig_2}
\end{figure}

The tested color videos are shown as Figure \ref{fig_1}, which are random selected from videoSegmentationData
dataset\footnote{\url{http://www.kecl.ntt.co.jp/people/kimura.akisato/saliency3.html}}. The last $32$ frames with size $288\times 352$ of each frame of the color videos are chosen to test. We study the best rank-$s$\footnote{Here and hereunder, rank-$s$ means TQt-rank-$s$ for TQt-SVD and Qt-rank-$s$ for Qt-SVD.} ($s=5,10,20$) approximations to the tested color videos and compare with the latest Qt-SVD \cite{DBLP:journals/appml/QinMZ22}. For our TQt-SVD, we use $4$ different orthogonal quaternion transformations $\dot{\mathbf{T}}$s: a random unitary quaternion matrix, the quaternion Discrete Fourier Transformation (QDFT), the quaternion Discrete Cosine Transformation (QDCT), and a Data-driven matrix obtained by taking the conjugate transpose of the left factor matrix of the mode-$3$ unfolding of the tested color videos. Table \ref{table1} gives the average  Peak Signal-to-Noise Ratio (PSNR) of all the frames of each tested color video. Figure \ref{fig_2} visually shows the results of one randomly selected frame.
\begin{table}[htbp]
	\caption{The averge PSNR of all the frames of each tested color video for different rank-$s$ approximations. (\textbf{Bold} fonts denote the best performance; \underline{underline} ones represent the second-best results).}
	\centering
	\resizebox{16.3cm}{1.5cm}{
\begin{tabular}{|c|l|ccc|ccc|ccc|ccc|}
	\hline
\multicolumn{2}{|c|} {Color videos:} & &AN119T   &  &  & DO01$\_$013&  &  &DO01$\_$030 &  &  &DO01$\_$014&  \\\toprule
	\hline
\multicolumn{2}{|c|} {Rank: $s=$}& 5 &10  &20  & 5 & 10 &20  & 5 & 10 & 20 & 5 & 10 & 20 \\
	\hline
	\multicolumn{2}{|c|} {Qt-SVD }&21.692 &25.502& 29.396  &\underline{20.321} &\underline{23.835} & \underline{27.869}&20.850  &24.022 &27.865 &24.487 &29.685  &34.293  \\
	\hline
	\multirow{4}{*}{TQt-SVD}&$\dot{\mathbf{T}}=$Random &21.375&25.056 &29.074&19.972&23.601&27.745 &20.682  &23.795  &27.597  &24.181  &29.608  &\underline{34.424}  \\\cline{2-14}
	&$\dot{\mathbf{T}}=$QDFT&21.692 &25.502&29.396 &\underline{20.321}  &\underline{23.835}&\underline{27.869}&20.850  &24.022 &27.865&24.487 &29.685  &34.293 \\\cline{2-14}
	&$\dot{\mathbf{T}}=$QDCT &\underline{21.723}  &\underline{25.570}  &\underline{29.514}  &20.306 &23.779 &27.807  &\underline{20.868}  &\underline{24.045}  &\underline{27.880}  &\underline{24.508}  &\underline{29.720}  &34.248  \\\cline{2-14}
	&$\dot{\mathbf{T}}=$Data-driven &\textbf{21.735} &\textbf{25.604 } &\textbf{29.551}  & \textbf{20.349}  & \textbf{23.872} &\textbf{27.984}  &\textbf{20.875}  &\textbf{24.060}  &\textbf{27.933}  & \textbf{24.558 } &\textbf{29.884}  &\textbf{34.541}  \\\cline{2-14}
	\hline
\end{tabular}}
\label{table1}
\end{table}

\renewcommand{\thesubfigure}{\roman{subfigure}} \makeatletter
\renewcommand{\@thesubfigure}{(\thesubfigure)\space}
\renewcommand{\p@subfigure}{\thefigure} \makeatother
From Table \ref{table1} and Figure \ref{fig_2}, we can see the advantages of TQt-SVD over Qt-SVD. The TQt-SVD is more flexible and can choose different transformations $\dot{\mathbf{T}}$s. The data-driven transformation has the best effect.
\section{Conclusions}
\label{sec:5}
This paper has introduced a flexible transform-based singular value decomposition for any $L$th-order ($L\geq 3$) quaternion tensors. Compared with the latest Qt-SVD, our proposed TQt-SVD can handle higher-order quaternion tensor (not necessary third-order) by using any invertible quaternion transform (not necessary fixed one).  In the experiments, we have verified the effectiveness of the proposed TQt-SVD in the application of the best TQt-rank-$s$ approximation for color videos. For some used transformations, \emph{e.g.}, $\dot{\mathbf{T}}=$Data-driven, the TQt-SVD has a better performance than the latest Qt-SVD. However, there may be a more suitable transformation, which requires our next study.

\section*{Acknowledgment}
This work was supported by The Science and Technology Development Fund, Macau SAR (File no. FDCT/085/2018/A2) and University of Macau (File no. MYRG2019-00039-FST).
\bibliographystyle{plain}
\bibliography{Myrefer}
\end{sloppypar}
\end{document}